\documentclass[12pt]{amsart}
\usepackage[cp850]{inputenc}
\usepackage{graphics}

\newtheorem{theorem}{Theorem}
\newtheorem{lemma}[theorem]{Lemma}

\newtheorem{corollary}[theorem]{Corollary}

\theoremstyle{definition}

\theoremstyle{remark}

\usepackage{amssymb,amsmath}

\def\n{\nabla}

\newcommand{\m}{\mbox{$M$}}

\newcommand{\s}{\mbox{$\Sigma$}}
\newcommand{\R}{\mbox{${\mathbb R}$}}
\newcommand{\N}{\mbox{$\m^2\times\R_1$}}
\newcommand{\g}[2]{\mbox{$\langle #1 ,#2 \rangle$}}
\newcommand{\fle}{\mbox{$\rightarrow$}}
\newcommand{\rf}[1]{\mbox{(\ref{#1})}}
\newcommand{\rl}[1]{{~\ref{#1}}}
\newcommand{\nablabar}{\mbox{$\overline{\nabla}$}}

\newcommand{\f}{\mbox{$f:\Sigma^2\fle\N$}}

\def\beq{\begin{equation}}
\def\eeq{\end{equation}}

\begin{document}

\title[A local estimate for maximal surfaces in Lorentzian product spaces]
{A local estimate for maximal surfaces in Lorentzian product spaces}

\author{Alma L. Albujer}
\address{Departamento de Matem\'{a}ticas, Universidad de Murcia, E-30100 Espinardo, Murcia, Spain}
\email{albujer@um.es}
\thanks{A.L. Albujer was supported by FPU Grant AP2004-4087 from Secretar\'{\i}a de Estado de Universidades e
Investigaci\'{o}n, MEC Spain.}

\author{Luis J. Al\'\i as}
\address{Departamento de Matem\'{a}ticas, Universidad de Murcia, E-30100 Espinardo, Murcia, Spain}
\email{ljalias@um.es}
\thanks{This work was partially supported by MEC project MTM2007-64504, and Fundaci\'{o}n S\'{e}neca
project 04540/GERM/06, Spain. This research is a result of the activity developed within the framework of the Programme in Support of Excellence
Groups of the Regi\'{o}n de Murcia, Spain, by Fundaci\'{o}n S\'{e}neca, Regional Agency for Science and Technology (Regional Plan for Science and
Technology 2007-2010).}



\dedicatory{Dedicated to Professor Manfredo P. do Carmo on the occasion of his 80th birthday}


\begin{abstract}
In this paper we introduce a local approach for the study of maximal surfaces immersed into a Lorentzian product space of the form $\m^2\times\R_1$,
where $\m^2$ is a connected Riemannian surface and $\m^2\times\R_1$ is endowed with the product Lorentzian metric. Specifically, we
establish a local integral inequality for the squared norm of the second fundamental form of the surface, which allows us to derive an alternative
proof of our Calabi-Bernstein theorem given in \cite{AA}.
\end{abstract}

\maketitle

\section{Introduction}
\label{s1} Maximal surfaces in 3-dimensional Lorentzian manifolds, that is, spacelike surfaces with zero mean curvature, have become a research field
of increasing interest in recent years, both from mathematical and physical points of view. In fact, one of the most relevant global results for
maximal surfaces in Lorentzian geometry is the well-known Calabi-Bernstein theorem, which states that the only complete maximal surfaces in the
3-dimensional Lorentz-Minkowski space $\R^3_1$ are the spacelike planes.

This result was firstly proved by Calabi \cite{Ca} and extended later to arbitrary dimension by Cheng and Yau \cite{CY}. After that, several
extensions and generalizations of the Calabi-Bernstein theorem have been given, and several alternatives proofs have been provided. In particular, in
\cite{AP2} the second author jointly with Palmer introduced a new approach to the Calabi-Bernstein theorem in the Lorentz-Minkowski space $\R_1^3$
based on a local integral inequality for the Gaussian curvature of a maximal surface in $\R^3_1$ which involved the local geometry of the surface and
the image of its Gauss map. As an application of it, they provided a new proof of the Calabi-Bernstein theorem in $\R^3_1$. In this paper, we
generalize this local approach to the case of maximal surfaces in a product space $M^2 \times \R$, where $M^2$ is a connected Riemannian surface and
$M^2 \times \R$ is endowed with the product Lorentzian metric
\[
\g{}{}=\pi_M^\ast(\g{}{}_M)-\pi_\mathbb{R}^\ast(dt^2).
\]
Here $\pi_M$ and $\pi_\mathbb{R}$ stand for the projections from $M^2\times\R$ onto each factor and $\g{}{}_M$ is the
Riemannian metric on $M$. For simplicity, we will simply write
\[
\g{}{}=\g{}{}_M-dt^2,
\]
and we will denote by \N\ the 3-dimensional Lorentzian product manifold obtained in that way. Specifically, we will prove the following extension of
\cite[Theorem 1]{AP2}.

\begin{theorem}
\label{local} Let $\m^2$ be an analytic Riemannian surface with non-negative Gaussian curvature, $K_M\geq 0$, and let \f\ be a maximal surface in \N.
Let $p$ be a point of \s\ and $R>0$ be a positive real number such that the geodesic disc of radius $R$ about $p$ satisfies $D(p,R)\subset\subset\s$.
Then for all $0<r<R$ it holds that \beq \label{ineq} 0\leq \int_{D(p,r)}\|A\|^2\mathrm{d}\Sigma\leq c_r\frac{L(r)}{r\log{(R/r)}}, \eeq where $L(r)$
denotes the length of the geodesic circle of radius $r$ about $p$, and
\[
c_r=\frac{\pi^2(1+\alpha_r^2)^2}{4\alpha_r\arctan{\alpha_r}}>0.
\]
Here $$\alpha_r=\sup_{D(p,r)}\cosh\theta\geq 1,$$ where $\theta$ denotes the hyperbolic angle between $N$ and $\partial_t$ along \s.
\end{theorem}
In particular, when \s\ is complete then the local integral inequality \rf{ineq} provides an alternative proof of the following parametric version of
the Calabi-Bernstein type result for complete maximal surfaces in Lorentzian product spaces given by the authors in \cite[Theorem 3.3]{AA}.
\begin{corollary}
\label{CBtheorem} Let $\m^2$ be a (necessarily complete) analytic Riemannian surface with non-negative Gaussian curvature, $K_M\geq 0$. Then any
complete maximal surface $\Sigma^2$ in \N\ is totally geodesic. In addition, if $K_M>0$ at some point on \m, then \s\ is a slice $\m\times\{ t_0\}$,
$t_0\in\R$.
\end{corollary}

As another application of Theorem\rl{local}, at points of a maximal surface where the second fundamental form does not vanish, we are able to
estimate the maximum possible geodesic radius in terms of a local positive constant.
\begin{corollary}
\label{proplocal} Let $\m^2$ be an analytic Riemannian surface with non-negative Gaussian curvature and let \f\ be a maximal surface in \N\ which is
not totally geodesic. Assume that $p\in\s$ is a point with $\|A\|(p)\neq 0$ and let $r>0$ be a positive real number such that
$D_r=D(p,r)\subset\subset\s$. Then
\[
R\leq r e^{C_r}
\]
for every $R>r$ with $D(p,R)\subset\subset\s$, where
\[
C_r=\frac{c_rL(r)}{r\int_{D_r}\|A\|^2}>0
\]
is a local positive constant depending only on the geometry of $f|_{D(p,r)}$.
\end{corollary}

A similar estimate for stable minimal surfaces in 3-dimensional Riemannian surfaces with non-negative Ricci curvature was given by Schoen in
\cite{Sc}. See also \cite{AP} for another similar estimate given by the second author and Palmer for the case of non-flat
spacelike surfaces with non-negative Gaussian curvature and zero mean curvature in a flat 4-dimensional Lorentzian
space.

\section{Preliminaries}
A smooth immersion \f\ of a connected surface $\Sigma^2$ is said to be a spacelike surface if the induced metric via $f$ is a Riemannian metric on
$\Sigma$, which as usual is also denoted by $\g{}{}$. Observe that
$$
\partial_t=(\partial/\partial_t)_{(x,t)}, \quad x\in\m, t\in\R,
$$
is a unitary timelike vector field globally defined on the ambient spacetime \N. This allows us to consider the unique unitary timelike normal field
$N$ globally defined on $\Sigma$ which is in the same time-orientation as $\partial_t$, so that
$$
\g{N}{\partial_t}\leq -1<0 \quad \mathrm{on} \quad \Sigma.
$$
We will refer to $N$ as the future-pointing Gauss map of $\Sigma$, and we will denote by $\Theta:\s\fle (-\infty,-1]$ the smooth function on \s\
given by $\Theta=\g{N}{\partial_t}$. Observe that the function $\Theta$ measures the hyperbolic angle $\theta$ between the timelike future-pointing
vector fields $N$ and $\partial_t$ along \s, since $\cosh{\theta}=-\Theta$.

Let $\nablabar$ and $\nabla$ denote the Levi-Civita connections in \N\ and $\Sigma$, respectively, and let $A:T\Sigma\fle T\Sigma$ stands for the
shape operator (or second fundamental form) of $\Sigma$ with respect to its future-pointing Gauss map $N$. It is well known that the Gauss and
Weingarten formulae for the spacelike surface \f\ are given by \beq \label{gaussfor} \nablabar_XY=\nabla_XY-\g{AX}{Y}N \eeq and \beq \label{wein}
AX=-\nablabar_XN, \eeq for any tangent vector fields $X,Y\in T\Sigma$. The mean curvature of a spacelike surface \f\ is defined by
$H=-(1/2)\mathrm{tr}{A}$, and \f\ is said to be a maximal surface when $H$ vanishes on \s.

The Gauss equation of a spacelike surface \s\ describes its Gaussian curvature $K$ in terms of the shape operator and the curvature of the ambient
space and it is given by \beq \label{gausseq} K=\overline{K}-\mathrm{det}A, \eeq where $\overline{K}$ denotes the sectional curvature in
$M^2\times\R_1$ of the plane tangent to \s. On the other hand, if $\overline{R}$ stands for the curvature tensor of the Lorentzian product \N, then
the Codazzi equation of \s\ describes the tangent component of $\overline{R}(X,Y)N$, for any tangent vector fields $X,Y \in T\s$, in terms of the
derivative of the shape operator. Specifically, it is given by \beq \label{codazzieq} (\overline{R}(X,Y)N)^\top=(\n_XA)Y-(\n_YA)X, \eeq where $\n_XA$
denotes the covariant derivative of $A$, that is,
\[
(\n_XA)Y=\n_X(AY)-A(\n_XY).
\]

In the particular case where \f\ is a maximal surface, it is not difficult to see that the Gauss \rf{gausseq} and Codazzi \rf{codazzieq} equations
for \s\ become \beq \label{gaussmaximal} K=\kappa_M \Theta^2+\frac{1}{2}\|A\|^2 \eeq and \beq \label{codazzibis}
(\n_XA)Y=(\n_YA)X+\kappa_M\Theta(\g{X}{\partial_t^\top}Y-\g{Y}{\partial_t^\top}X), \eeq for any tangent vector fields $X,Y \in T\Sigma$,
respectively. Here $\|A\|^2=\mathrm{tr}(A^2)$ and $\kappa_M$ stands for the Gaussian curvature of $M$ along the surface \s, that is, $\kappa_M=K_M
\circ \Pi \in \mathcal{C}^\infty(\s)$ where $K_M$ is the Gaussian curvature of $M$ and $\Pi=\pi_M \circ f:\s \fle M$ denotes the projection of \s\
onto $M$. Here and in what follows, $Z^\top\in T\s$ denotes the tangential component of a vector field $Z$ along the immersion \f, that is
\[
Z=Z^\top-\g{N}{Z}N.
\]
Thus, in particular, \beq \label{desct}\partial_t^\top=\partial_t+\Theta N, \eeq (for the details see \cite{AA}). Taking norms in the last expression
we get \beq \label{eq6} \|\partial_t^\top\|^2=\Theta^2-1. \eeq

It is well known that a spacelike surface \f\ is locally a spacelike graph over $M$ (see for instance \cite[Lemma 3.1]{AA}), that is, for any given
point $p \in \s$, there exists an open subset $\Omega$ on $M$ containing $\Pi(p)$, $\Pi(p) \in \Omega\subset M$, and a function $u \in
\mathcal{C}^\infty(\Omega)$ such that the surface \s\ is locally given in a neighborhood of $p$ by $\s(u)=\{(x,u(x)):x\in\Omega\}\subset \N$.
Therefore, the metric induced on $\Sigma(u)$ from the Lorentzian metric on the ambient space is given by
\beq
\label{gu}
\g{}{}=\g{}{}_M-du^2.
\eeq
The condition that $\Sigma(u)$ is spacelike becomes $|Du|^2<1$ on $\Omega\subset M$, where $Du$ denotes the gradient of $u$ in $M$
and $|Du|$ denotes its norm. Finally, it is not difficult to see that the mean curvature function $H$ of $\Sigma(u)$
is given by
\[
2H=\mathrm{Div}\left( \frac{Du}{\sqrt{1-|Du|^2}}\right),
\]
on $\Omega$, where $\mathrm{Div}$ stands for the divergence operator on $M$ with respect to the metric $\g{}{}_M$. In
particular, a spacelike immersion \f\ is a maximal surface if and
only if it is locally given as the graph of a function $u$ satisfying the following partial differential equation,
\beq
\label{zmc}
\mathrm{Div}\left(\frac{Du}{\sqrt{1-|Du|^2}}\right)=0, \quad |Du|^2<1.
\eeq

\section{Proof of the results}

The proof of Theorem\rl{local} is inspired by the ideas in \cite{AP2}, and it is an application of the following intrinsic property.
\begin{lemma}\cite[Lemma 3]{AP2}
\label{lemma} Let \s\ be an analytic Riemannian surface with non-negative Gaussian curvature $K\geq 0$. Let $\psi$ be a smooth function on \s\ which
satisfies
\[
\psi\Delta\psi\geq 0
\]
on \s. Then for $0<r<R$
\[
\int_{D_r}\psi\Delta\psi\leq\frac{2 L(r)}{r\log{(R/r)}} \sup_{D_R}\psi^2,
\]
where $D_r$ denotes the geodesic disc of radius $r$ about a fixed point in \s, $D_r\subset D_R\subset\subset\s$, and $L(r)$ denotes the length of
$\partial D_r$, the geodesic circle of radius r.
\end{lemma}

\begin{proof}[Proof of Theorem\rl{local}]
Observe that since \m\ is analytic and \s\ is locally given by the maximal surface equation \rf{zmc}, then \s, endowed with the induced metric, is
also an analytic Riemannian surface. Besides, from \rf{gaussmaximal} we also know that the Gaussian curvature of \s\ is non-negative, $K\geq 0$.
Therefore, we may apply Lemma\rl{lemma} to an appropriate smooth function $\psi$. Let us consider $\psi=\arctan\Theta$.

Since $\partial_t$ is parallel on \N\ we have that \beq \label{eq4} \nablabar_X\partial_t=0 \eeq for any tangent vector field $X\in T\Sigma$. Thus,
\[
X(\Theta)=\g{\nablabar_XN}{\partial_t}=-\g{AX}{\partial_t^\top}=-\g{X}{A\partial_t^\top}
\]
for every $X\in T\Sigma$, and then the gradient of $\Theta$ on \s\ is given by \beq \label{gradT} \nabla\Theta=-A\partial_t^\top. \eeq Therefore,
from \rf{gradT} and \rf{eq6} we obtain \beq \label{normamaximal} \|\nabla\Theta\|^2=\frac{1}{2}\|A\|^2(\Theta^2-1), \eeq since for a maximal surface
it holds $A^2=(1/2)\|A\|^2I$.

On the other hand, taking into account \rf{desct}, and using Gauss \rf{gaussfor} and Weingarten \rf{wein} formulae, \rf{eq4} also yields \beq
\label{eq5} \n_X\partial_t^\top=-\Theta AX \eeq for every $X \in T\s$. Therefore, using Codazzi equation \rf{codazzibis} and equations \rf{eq6} and
\rf{eq5} we get
\begin{eqnarray*}
\nabla_X\nabla\Theta & = & -(\nabla_XA)(\partial_t^\top)-A(\nabla_X\partial_t^\top) \\
{} & = & -(\nabla_{\partial_t^\top}A)(X)-\kappa_M\Theta\left(\g{X}{\partial_t^\top}\partial_t^\top-\|\partial_t^\top\|^2X\right)
+\Theta A^2X \\
{} & = & -(\nabla_{\partial_t^\top}A)(X)+\kappa_M\Theta\left((\Theta^2-1)X-\g{X}{\partial_t^\top}\partial_t^\top\right) +\Theta A^2X,
\end{eqnarray*}
for every $X\in T\Sigma$. Thus, the Laplacian of $\Theta$ is given by \beq \label{laplaT} \Delta\Theta = \Theta(\kappa_M(\Theta^2-1)+\|A\|^2), \eeq
since
\[
\mathrm{tr}(\nabla_{\partial_t^\top}A)=\nabla_{\partial_t^\top}(\mathrm{tr}A)=0.
\]

Using \rf{laplaT} and \rf{normamaximal} we can compute
\[
\Delta\psi=\frac{\Delta\Theta}{1+\Theta^2}-\frac{2\Theta\|\nabla\Theta\|^2}{(1+\Theta^2)^2}=
\frac{2\Theta}{(1+\Theta^2)^2}\|A\|^2+\frac{(\Theta^2-1)\Theta}{1+\Theta^2}\kappa_M,
\]
and therefore, taking into account that $\Theta\arctan\Theta\geq 0$, $\Theta\leq -1$ and $\kappa_M\geq 0$, we obtain \beq \label{eq17}
\psi\Delta\psi=\frac{2\Theta\arctan\Theta}{(1+\Theta^2)^2}\|A\|^2+ \frac{(\Theta^2-1)\Theta\arctan\Theta}{1+\Theta^2}\kappa_M\geq
\phi(\Theta)\|A\|^2, \eeq where
\[
\phi(s)=\frac{2s\arctan s}{(1+s^2)^2}.
\]
Observe that the function $\phi(s)$ is strictly increasing for $s\leq -1$. Since $-\alpha_r\leq\Theta\leq-1$ on $D(p,r)$, we get
\[
\phi(\Theta)\geq\phi(-\alpha_r)=\frac{2\alpha_r\arctan\alpha_r}{(1+\alpha_r^2)^2} \quad \mathrm{on} \quad D(p,r),
\]
which, jointly with \rf{eq17}, yields
$$
\psi\Delta\psi\geq\frac{2\alpha_r\arctan\alpha_r}{(1+\alpha_r^2)^2}\|A\|^2  \quad \mathrm{on} \quad D(p,r).
$$
Integrating now this inequality over $D(p,r)$ and using Lemma\rl{lemma} we conclude that
\[
0\leq\frac{2\alpha_r\arctan\alpha_r}{(1+\alpha_r^2)^2}\int_{D(p,r)}\|A\|^2\mathrm{d}\Sigma\leq \int_{D(p,r)}\psi\Delta\psi\leq
\frac{\pi^2}{2}\frac{L(r)}{r\log{(R/r)}},
\]
which yields \rf{ineq}.
\end{proof}

\begin{proof}[Proof of Corollary\rl{CBtheorem}]
Since \s\ is complete, then $R$ can approach to infinity in \rf{ineq} for a fixed arbitrary $p\in\s$ and a fixed $r$, which gives
\[
\int_{D(p,r)}\|A\|^2\mathrm{d}\Sigma=0.
\]
Therefore, $\|A\|^2=0$ and \s\ must be totally geodesic. From \rf{gradT}, this implies that $\Theta=\Theta_0\leq-1$ is constant on \s, and then
\rf{laplaT} implies that, when $K_M>0$ somewhere in \m, it must be $\Theta_0=-1$. Finally, by \rf{eq6} we conclude that \s\ must be a slice.
\end{proof}

Corollary\rl{proplocal} is a direct consequence of Theorem\rl{local}.

\bibliographystyle{amsplain}

\end{document}